\documentclass[11pt,twoside,english]{myamsart}
\normalsize
\advance\oddsidemargin by -1.0cm
\advance\evensidemargin by -1.0cm
\advance\topmargin by -1.0cm 
\usepackage{amssymb}
\usepackage{babel}
\usepackage{amsmath}
\usepackage{amscd}   
\usepackage{epsfig}  
\usepackage{rotating}
\usepackage{psfrag}
\theoremstyle{plain}
\newtheorem{cor}{Corollary}[section]

\newtheorem{thm}{Theorem}[section]
\newtheorem{prop}{Proposition}[section]
\theoremstyle{definition}
\newtheorem{exa}{Example}[section]
\newtheorem{NB}{Remark}[section]

\renewcommand{\labelenumi}{$(\theenumi)$}
\newcommand{\bdm}{\begin{displaymath}}
\newcommand{\edm}{\end{displaymath}}
\newcommand{\be}{\begin{equation}}
\newcommand{\ee}{\end{equation}}
\newcommand{\ba}[1]{\begin{array}{#1}}
\newcommand{\ea}{\end{array}}
\newcommand{\bea}[1][]{\begin{eqnarray#1}}
\newcommand{\eea}[1][]{\end{eqnarray#1}}
\newcommand{\btab}{\begin{tabular}}
\newcommand{\etab}{\end{tabular}}

\newcommand{\op}{\oplus}

\newcommand{\ra}{\rightarrow}

\newcommand{\lan}{\left\langle}
\newcommand{\ran}{\right\rangle}

\newcommand{\vol}{\ensuremath{\mathrm{vol}}}

\newcommand{\Id}{\ensuremath{\mathrm{Id}}}

\newcommand{\R}{\ensuremath{\mathbb{R}}}

\newcommand{\vphi}{\ensuremath{\varphi}}    

\newcommand{\hut}{\wedge}

\newcommand{\Ric}{\ensuremath{\mathrm{Ric}}}
\newcommand{\Scal}{\ensuremath{\mathrm{Scal}}}
\newcommand{\Scalg}{\ensuremath{\mathrm{Scal}^g}}

\begin{document}
\def\haken{\mathbin{\hbox to 6pt{%
                 \vrule height0.4pt width5pt depth0pt
                 \kern-.4pt
                 \vrule height6pt width0.4pt depth0pt\hss}}}
    \let \hook\intprod
\setcounter{equation}{0}
%
%------ draw title page -----
%
\thispagestyle{empty}
%
%\hbox to \hsize{%
%  \vtop{} \hfill
%  \vtop{\hbox{PRELIMINARY VERSION}}}
%------------------------------
\date{\today}
%---------------------------------------------------------------------------
\title[Eigenvalue estimates for Dirac operators with parallel 
torsion]{Eigenvalue estimates for Dirac operators with parallel
characteristic torsion}
%---------------------------------------------------------------------------
%
% author and address
%
%-------------------------------------------
%
\author{Ilka Agricola}
\author{Thomas Friedrich}
\author{Mario Kassuba}
%-------------------------------------------
\address{\hspace{-5mm} 
{\normalfont\ttfamily agricola@mathematik.hu-berlin.de}\newline
{\normalfont\ttfamily friedric@mathematik.hu-berlin.de}\newline
{\normalfont\ttfamily kassuba@mathematik.hu-berlin.de}\newline
Institut f\"ur Mathematik \newline
Humboldt-Universit\"at zu Berlin\newline
Unter den Linden 6\newline
Sitz: John-von-Neumann-Haus, Adlershof\newline
D-10099 Berlin, Germany}
% 
%-----------------------------------------------------------
\thanks{Supported by the SFB 647 "Space---Time---Matter" of the DFG and
the Junior Research Group "Special Geometries in Mathematical Physics"
of the Volkswagen\textbf{Foundation}.}
%-----------------------------------------------------------
\subjclass[2000]{Primary 53 C 25; Secondary 81 T 30}
%-----------------------------------------------------------
%-----------------------------------------------------------
\begin{abstract}
%---------------
Assume that the compact Riemannian spin manifold $(M^n,g)$ admits a 
$G$-structure with characteristic connection $\nabla$ and
 parallel characteristic 
torsion ($\nabla T=0$), and consider the Dirac operator $D^{1/3}$ 
corresponding to the torsion $T/3$. This operator  plays an eminent 
role in the investigation of such manifolds and includes as special
cases Kostant's ``cubic Dirac operator'' and the Dolbeault operator.
In this article, we describe a general method of computation for lower 
bounds of the eigenvalues of $D^{1/3}$ by a clever deformation of
the spinorial connection.   In order to get explicit bounds, each 
geometric structure needs to be investigated separately; we
do this in full generality in dimension $4$ and for Sasaki manifolds
in dimension $5$.
\end{abstract}
%-------------
\maketitle
%----------------
%\tableofcontents
%----------------
\pagestyle{headings}
%
%
%-------------- body of the document ------------------------------------------
%
%---------------------------------------------------------------------------- 
\section{Introduction}\noindent
%----------------------------------------------------------------------------
Lower bounds for the first eigenvalue of the Riemannian Dirac operator $D^g$ on
a compact Riemannian spin manifold depending on the scalar curvature
are well known since more than two decades (see \cite{Friedrich80}, 
\cite{Kirchberg}). In  past years, another operator of Dirac type turned out 
to play a crucial role  in the investigation
of non-integrable geometric structures as well as in several models in
superstring theory. Indeed, for many of these geometries it is known that
one can replace the Levi-Civita connection by a unique adapted metric 
connection $\nabla$ with skew-symmetric torsion $T$ preserving the
geometric structure, the so-called \emph{characteristic connection} of
the geometric structure. The survey article \cite{Agricola06} discusses these
developments. The Dirac operator in question is then not merely the
Dirac operator associated with $\nabla$, but the operator
$D^{1/3} = D^g + \frac{1}{4}T$ corresponding to the torsion form
$T/3$ (see  \cite{Agricola&F04a}, \cite{Agricola&F04b}).
In fact, $D^{1/3}$ coincides with the so-called
``cubic Dirac operator'' introduced by B.~Kostant (\cite{Kostant99}, 
\cite{Agricola03}) on naturally reductive spaces and with the
Dolbeault operator on Hermitian manifolds (\cite{Bismut}, \cite{Gauduchon97}).

The aim of the present paper is to estimate the eigenvalues of this new 
Dirac operator in case the torsion form $T$ is $\nabla$-parallel. In this
situation $(D^{1/3})^2$ and the torsion form $T$ commute on spinors (see 
\cite{Agricola&F04b}) and we can estimate the eigenvalues separately in any
eigensubbundle of the symmetric endomorphism defined by $T$. The classically
well-known classes of manifolds with parallel characteristic torsion---nearly 
K\"ahler manifolds, Sasakian manifolds, nearly parallel $G_2$-manifolds, 
naturally reductive spaces---have been considerably enlarged in more recent
investigations (see \cite{Vaisman79}, \cite{Gauduchon&O98}, 
\cite{Friedrich&I1}, \cite{Alexandrov03}, \cite{Alexandrov&F&S04}, 
\cite{Agricola&F06}, 
\cite{Friedrich06}, \cite{Schoemann}), leading eventually to an abundant
supply of manifolds to which our results can be applied.
It was known that the general formula of 
Schr\"odinger-Lichnerowicz type (S-L-formula for short) for
the operator $D^{1/3}$ derived in \cite{Agricola&F04a}, \cite{Agricola&F04b}
does not yield optimal lower bounds for the spectrum. In this article,
we deform the connection $\nabla$ by polynomials of the torsion form.
The resulting connections are not affine connections anymore, they
exist only on the spinor bundle. In 1980, Friedrich had used a similar
spinorial  modification of the lift of the Levi-Civita connection to
derive his estimate for the eigenvalue of $D^g$; the main difference 
is, however, that there was no torsion form to use there.

The article is organized as follows. In section $2$, we
prove the necessary integral formulas for perturbations of $D^{1/3}$ 
by some parallel symmetric endomorphism $S$ of the spinor bundle and we 
describe the  general strategy for proving bounds of the spectrum 
of $D^{1/3}$. 
In order to obtain an explicit estimate through this method, one
needs to know the algebraic type of $T$ and the splitting of the spinor 
bundle, hence, every special geometry requires a separate investigation. 
Section $3$ is devoted to the determination of the lower bound for
the first eigenvalue $\lambda$ of $(D^{1/3})^2$ by this method on 
$4$-dimensional compact spin manifolds with positive scalar curvature
$\Scalg_{\min}>0$ and 
parallel torsion $T\neq 0$. This applies for example to
generalized  Hopf manifolds, i.\,e.~Hermitian $4$-manifolds with Lee form
parallel with respect to the Levi-Civita connection (see \cite{Vaisman79}, 
\cite{Gauduchon&O98},  \cite{Belgun00}). 
We prove  the 
following estimate depending on the ratio $c:=\Scalg_{\min} / ||T||^2>0$:
\bdm
\lambda \ \geq\ \left\{\ba{ll} 
\frac{1}{16}\left[\sqrt{6\,\Scal^g_{\min}} - ||T||\right]^2\, = \,
\frac{||T||^2}{16}\left[\sqrt{6\,c}-1\right]^2 & \text{for } 
1/6\leq c\leq 3/2   
\\[2mm]
\frac{1}{4}\left[\Scalg_{\min}-\frac{1}{2}||T||^2\right] \, = \, 
\frac{||T||^2}{4}\left[c -\frac{1}{2}\right]
& \text{for } c\geq 3/2 
\ea\right.
\edm
We refer to Theorem \ref{thm-n4} for details.
 In section $4$ we apply our method 
to $5$-dimensional compact Sasakian manifolds with $\Scalg_{\min}> -4$. 
Their characteristic torsion 
is always parallel of length $8$ and we get the optimal estimate
\bdm
\lambda \ \geq\ \left\{\ba{ll} 
\frac{1}{16}\left[1+\frac{1}{4}\Scalg_{\min}\right]^2 & \text{ for } 
-4<\Scalg_{\min}\leq 4(9+4\sqrt{5})\\[2mm]
\frac{5}{16}\, \Scalg_{\min} & \text{ for }\ \Scalg_{\min}\geq 4(9+4\sqrt{5})
\simeq 71,78.
\ea\ \right.
\edm
Full statements are to be found in Theorem \ref{est-dim5}.
The lower bound is attained, for example, on every $\eta$-Einstein-Sasakian 
manifold.  It is a curious fact that this is, to our knowledge, the first
eigenvalue estimate for a Dirac operator with a \emph{quadratic} dependence
on the scalar curvature. In both dimensions, we have the effect that
the main improvement is on an intervall corresponding to `small' scalar
curvatures and eigenvalues (upper line in the estimates above).
Beyond this bound, the estimate in dimension $4$ can be obtained relatively 
easy by universal arguments. In dimension $5$, this bound is exactly
the Riemannian estimate $n/4(n-1)\Scalg_{\min}$, a
non-trivial fact that can be traced back to the particular property that 
$0$ is an eigenvalue of the characteristic torsion $T$ of a Sasaki
manifold (this can never happen in dimension $4$).
%
%---------------------------------------------------------------------------- 
\section{Schr\"odinger-Lichnerowicz formulas for the deformation of the 
connection}\label{fam-conn}\noindent
%----------------------------------------------------------------------------
%
Consider a compact Riemannian spin manifold $(M^n,g, T)$ with Levi-Civita 
connection
$\nabla^g$ as well as a  metric connection with skew-symmetric 
torsion $T\in\Lambda^3(M)$,
\bdm
\nabla_X Y\ :=\ \nabla^{g}_X Y + \frac{1}{2} \cdot T(X,Y,-)\,.
\edm
This connection can  be lifted to the spinor
bundle $\Sigma M$ of $M$, where it takes the expression
\bdm
\nabla_X \psi\ :=\ \nabla^{g}_X \psi + \frac{1}{4} ( X\haken T)\cdot\psi\,.
\edm
Its Dirac operator is given by $D=D^g+ (3/4)\,T$, where $D^g$ denotes 
the Riemannian Dirac operator. Besides this, the connection with
torsion $T/3$---henceforth denoted $\nabla^{1/3}$---will play a crucial
role in our considerations. Its associated Dirac operator is accordingly
given by $D^{1/3}=D^g + T/4= D-T/2$. Similarly, the spinor Laplacian of
$\nabla$ will be written $\Delta$. The Laplacian of $\nabla^{1/3}$
will never be used. We define an algebraic $4$-form
derived from $T$ by 
\bdm
\sigma_{T} \ :=\ \frac{1}{2} \sum_k (e_k\haken T)\hut
(e_k\haken T) \, ,
\edm
where $e_1,\ldots,e_n$ is an orthonormal basis. Recall that any $k$-form
acts on spinors by the extension of Clifford multiplication; henceforth,
we shall make no notational difference between a $k$-form and the 
endomorphism on the spinor bundle that it induces.
\begin{prop}[{\cite{Agricola&F04a}}]\label{3-form-square}
%-------------------------------------------------------
Let $T$ be a $3$-form in dimension $n\geq 5$, and denote by the same symbol 
its associated $(2,1)$-tensor. Then its square inside
the Clifford algebra has no contributions of degree $2$ and $6$, and its scalar
and fourth degree part are given by
\bdm
T^2_0\ =\ \frac{1}{6}\,\sum_{i,j=1}^n ||T(e_i,e_j)||^2\ =:\ || T ||^2,\quad
T^2_4\ =\ - \, 2 \, \sigma_T.
\edm
For $n=3,4$, one has $T^2=||T||^2$.
\end{prop}
\noindent
Let's now state the main S-L-formula for $(D^{1/3})^2$.
It links the Dirac operator for the torsion $T/3$ with the Laplacian 
for the torsion $T$. Here, $\Scalg$ and $\Scal$ denote the scalar
curvatures of the Levi-Civita connection and the new connection $\nabla$,
respectively. They are related by $\Scal=\Scalg-(3/2)||T||^2$.
\begin{thm}[{\cite{Agricola&F04a}}]\label{new-weitzenboeck}
%---------------------------------------------------------
The spinor Laplacian $\Delta$ and the square of the Dirac operator
$D^{1/3}$ are related by
\bdm
(D^{1/3})^2\ =\ \Delta + \frac{1}{4} \, dT 
+\frac{1}{4}\,\Scalg - \frac{1}{8}\, || T||^2.
\edm
\end{thm}
\noindent
Eigenvalue estimates for $D^{1/3}$ footing on this relation will
be called \emph{universal}, to distinguish them
from the new estimates obtained by deforming $\nabla$ to be discussed
later.

\noindent
In this article, we are only interested in connections $\nabla$ for
which $T$ is parallel, $\nabla T=0$. In this case, $T$ has constant length
 and it is well-known
that $dT=2 \, \sigma_T$, hence Proposition~\ref{3-form-square}
implies
\be\label{dT}
dT \ =\ - T^2 + ||T||^2.
\ee 
Combined with the main result of Theorem \ref{new-weitzenboeck}, we obtain
in the case of parallel torsion the ``universal'' S-L-formula
\be\label{D/3-square}
(D^{1/3})^2\ =\ \Delta -\frac{1}{4}\,T^2 
+\frac{1}{4}\,\Scalg + \frac{1}{8} || T||^2.
\ee
Now let $S:\Sigma M\ra\Sigma M$ be a symmetric endomorphism that is  
parallel itself, $\nabla S=0$. The main case we have in mind in our 
applications are polynomials $S = P(T)$ in $T$.
Then we can define a new, \emph{$S$-deformed spin connection} $\nabla^S$
on  $\Sigma M$ by
\bdm
\nabla^S_X \psi \ :=\ \nabla_X\psi - \frac{1}{2}(X\cdot S+ S\cdot X)\cdot\psi
\edm
which is metric again. Indeed, the symmetry  of $S$ implies for
any two spinors $\vphi,\psi\in\Sigma M$
\be\label{antisym}
\langle(X\cdot S+ S\cdot X)\,\vphi,\psi \rangle +
\langle\vphi, (X\cdot S+ S\cdot X)\,\psi \rangle\ =\ 0.
\ee
In the following technical proposition, we gather all necessary computations
involving $S$. Here, all lengths and inner products refer to  the
$L^2$-scalar product on spinors.
\begin{prop}\label{general-formulas}
%-----------------------------------
\begin{enumerate}
\item[]
\item[a)] $\displaystyle D^{1/3} S+ SD^{1/3}=\sum_{i=1}^n (e_i\cdot S+S\cdot
  e_i)\nabla_{e_i}-\frac{1}{2}(TS+ST)$
\item[b)] $\displaystyle ||\nabla^S\psi||^2=||\nabla\psi||^2 +\sum_{i=1}^n \langle 
(e_i\cdot S+S\cdot e_i)\nabla_{e_i}\psi,\psi\rangle+\frac{1}{4}\sum_{i=1}^n
||(e_i\cdot S+S\cdot e_i)\psi ||^2$
\item[c)] $\displaystyle\langle (D^{1/3}+S)^2\psi,\psi\rangle =
||\nabla^S\psi||^2-\frac{1}{4}\sum_{i=1}^n||(e_i\cdot S+S\cdot e_i)\psi ||^2-
\frac{1}{4}|| T\psi||^2+\frac{1}{8}||T||^2\cdot ||\psi||^2 +
\frac{1}{4} \int_{M^n}\Scalg\,||\psi||^2 + ||S\psi||^2-
\langle TS\,\psi,\psi\rangle $ 
\end{enumerate}
\end{prop}
\begin{proof}
%-------------
Identity a) is easy:
\bea[*]
\left[D^{1/3} S+ SD^{1/3}\right]\psi &= & \left[(D- T/2)S+S(D- T/2)\right]\psi
\, =\, DS\psi+SD\psi-\frac{1}{2}(TS+ST)\psi\\ 
&=& \sum_{i=1}^n e_i\cdot\nabla_{e_i} (S\psi)+ 
\sum_{i=1}^n S e_i\nabla_{e_i}\psi -\frac{1}{2}(TS+ST)\psi.
\eea[*]
Since $S$ is assumed to be parallel, the claim follows. The second 
relation is a direct consequence of the definition of $\nabla^S$ and
the antisymmetry property stated in  equation $(\ref{antisym})$.

For the last statement, observe that $D^{1/3}+S$ is again a symmetric
first order differential operator. First, we have
\bdm
(D^{1/3}+S)^2 \ = \ (D^{1/3})^2+ (D^{1/3}S+SD^{1/3})+S^2 .
\edm
Hence, we can insert identity $(\ref{D/3-square})$ and relation a)
\bdm
(D^{1/3}+S)^2 \ = \ \Delta -\frac{1}{4}\,T^2 
+\frac{1}{4}\,\Scalg + \frac{1}{8} || T||^2 + \sum_{i=1}^n (e_i\cdot S+S\cdot
  e_i)\nabla_{e_i}-\frac{1}{2}(TS+ST) +S^2,
\edm
which implies for the scalar product
\bea[*]
\langle (D^{1/3}+S)^2\psi,\psi\rangle& =&
\sum_{i=1}^n \langle (e_i\cdot S+S\cdot e_i)\nabla_{e_i}
\psi,\psi\rangle
-\frac{1}{4}\,||T\psi||^2 +\frac{1}{4}\int_{M^n}\Scalg \,||\psi||^2\\
& & + \frac{1}{8} || T||^2\cdot ||\psi||^2 + ||S\psi||^2 - \frac{1}{2}\langle
(TS+ST)\psi,\psi\rangle + ||\nabla \psi||^2     .
\eea[*]
Now the result follows from relation b). For the very last term,
observe that 
\bdm
\langle TS\psi,\psi\rangle\,=\,\langle S\psi,T\psi\rangle\, =\, 
\langle \psi,ST\psi\rangle
\edm
and that the imaginary part is irrelevant when integrating.
\end{proof}
\begin{NB}
%---------
We will refer to formula c) as the 
\emph{$S$-deformed  Schr\"odinger-Lichnerowicz formula}.
\end{NB}
\noindent
The general strategy for deducing  eigenvalue estimates for any eigenvalue
$\lambda$ of $(D^{1/3})^2$ from the deformed
connection $\nabla^S$ is as follows. The torsion $T$ is a 
$\nabla$-parallel symmetric endomorphism of the spinor bundle. We split
the spinor bundle under $T$,
\bdm
\Sigma M \ = \ \bigoplus_{\mu} \Sigma_{\mu} \, .
\edm
The connection $\nabla$ preserves this splitting. Moreover, the fact that
 $T$ and $(D^{1/3})^2$ commute (see \cite[Prop.~3.4]{Agricola&F04b})
implies that $(D^{1/3})^2$ preserves this splitting, too. Remark that,
in general, the first order operator $D^{1/3}$ does not commute with $T$.
Next we make
a suitable Ansatz for the endomorphism $S$. A natural restriction on
$S$ is that the new connection $\nabla^S$ should  preserve
the decomposition of $\Sigma M$ again. 
If $S$ is a polynomial in $T$, we obtain 
a simple characterization of such endomorphisms. Indeed, let us fix one 
subbundle $\Sigma_{\mu_0}$ and consider the minimal number of indices 
$\mu_1, \ldots , \mu_k$
such that
\bdm
TM^n \cdot \Sigma_{\mu_0} \ = \ \big\{ X \cdot \psi \, : \, X \in TM^n \, ,
\, \psi \in  \Sigma_{\mu_0} \big\} \ \subset \ \Sigma_{\mu_1} \oplus
\ldots \oplus \Sigma_{\mu_k} \, .
\edm
Then an easy computation yields the following result.
\begin{prop} \label{polynom}
%-----------
Let $S = P(T)$ be a polynomial in the torsion form $T$. Then $\nabla^S$ is a
connection in the subbundle $\Sigma_{\mu_0}$ if and only if, for
all $1\leq i\leq k$:
\bdm
\big(P(\mu_0) \, + \ P(\mu_i) \big) \, \mu_i \ =  \ 
\big(P(\mu_0) \, + \ P(\mu_i) \big) \, \mu_0  .  
\edm
\end{prop}
\noindent
Maximizing over all admissible endomorphisms $P(T)=S $ in a fixed subbundle
$\Sigma_{\mu_0}$, the $S$-deformed
S-L-formula  yields a lower
bound for the first eigenvalue $\lambda$ of $(D^{1/3})^2$
on $\Sigma_{\mu_0}$. However,  this estimate typically still involves the
eigenspinor $\psi$ of  $(D^{1/3})^2$, in particular the term 
$\langle D^{1/3}\psi,\psi\rangle$. In general, 
$\psi \in \Gamma(\Sigma_{\mu_0})$ will not be an eigenspinor of $D^{1/3}$,
hence this term cannot be expressed through $\lambda$ in any simple way. 
Nevertheless, it will be possible to overcome this difficulty in
special geometric situations by some additional arguments.
Lower estimates for $\lambda$ obtained by this strategy will,
for better reference, be called \emph{$S$-deformed eigenvalue bounds}.
In the next sections, we will show that this strategy yields 
non-trivial new eigenvalue estimates in concrete special geometries.
%
%-------------------------------------------------------------------------
\section{$4$-dimensional manifolds with  parallel torsion}
%-------------------------------------------------------------------------
%
\noindent
A particular property of the dimension $4$ is that $\sigma_T=0$ for purely 
algebraic reasons, hence $\nabla T=0$ implies $d T=0$ and $T^2$ acts
by scalar multiplication with $||T||^2$. 
Thus equation  $(\ref{D/3-square})$
reduces to
\bdm
(D^{1/3})^2\ =\ \Delta +\frac{1}{4}\,\Scalg - \frac{1}{8} || T||^2,
\edm
and we obtain for any eigenvalue $\lambda$ of $(D^{1/3})^2$
the following universal lower bound 
\be\label{dim4}
\lambda\ \geq\ \frac{1}{4}\left[\Scalg_{\min}-\frac{1}{2}||T||^2\right].
\ee
Remark that even for $T=0$ a better estimate for the first eigenvalue
of the Riemannian Dirac operator is known, $\lambda \geq \Scalg_{\min}/3$, see 
\cite{Friedrich80}.
Equality is obtained if and only
if $\nabla$ admits at least one parallel spinor $\psi\neq 0$. A look at
the  usual integrability condition $\nabla\nabla\psi=0$
and $\sigma_T=\nabla T=0$ yields the $\Ric^\nabla$-flatness 
of the connection $\nabla$. In particular, $0=\Scal=\Scalg -3||T||^2/2$
holds in this situation, while in general,
$\Scalg_{\min}$ and $||T||$ are independent geometric quantities.
Obviously, the universal bound becomes useless for 
$2\,\Scalg_{\min}\leq ||T||^2$. For complex  Hermitian spin surfaces of
non-negative conformal scalar curvature, a discussion of the universal bound
may be found in \cite[Thm 1.2]{Alexandrov&G&I01}.
\\

\begin{figure}
\bdm
\psfrag{1/6}{$\frac{1}{6}$}\psfrag{1/2}{$\frac{1}{2}$}
\psfrag{3/2}{$\frac{3}{2}$}\psfrag{1}{$1$}
\psfrag{2}{$2$}\psfrag{t2/5}{$\frac{1}{5}||T||^2$}
\psfrag{2t2/5}{$\frac{2}{5}||T||^2$}\psfrag{c}{$c$}
\includegraphics[width=7cm]{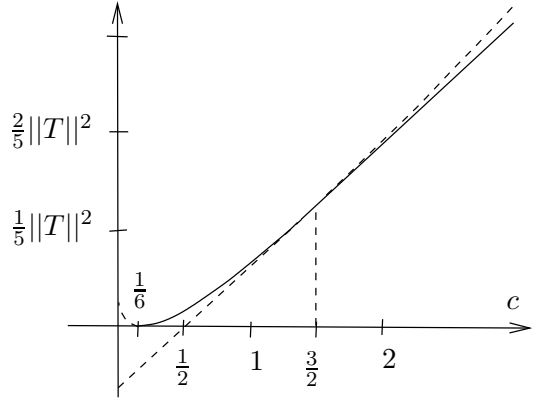}
\edm
\caption{Lower bound for $\lambda$ according to Theorem \ref{thm-n4}.}
\label{new-estimate}
\end{figure}

\noindent
Consider the symmetric endomorphism 
$S:=a_0\,\Id +a_1\, T:\Sigma M\ra\Sigma M$ for two real constants $a_0, a_1$, 
which is symmetric and again parallel. Higher order terms in $T$ are not 
needed, as $T^2$ is already acting by a scalar. In particular, the
spinor bundle splits into two subbundles of equal dimension, 
\bdm
\Sigma M \ =\ \Sigma_-\oplus\Sigma_+, \quad\Sigma_{\pm}\ :=\
\{\psi\in \Sigma M\, :\, T\cdot\psi = \pm ||T||\,\psi \}.
\edm
We shall use the strategy outlined in Section \ref{fam-conn} to prove the
following result.
\begin{thm}\label{thm-n4}
%------------------------
Let $(M^4,g)$ be a compact, $4$-dimensional spin manifold with
$\Scalg_{\min}>0$ and 
$0\neq T\in \Lambda^3(M)$ a $3$-form such that the connection $\nabla$ it
defines satisfies $\nabla T=0$.  The first eigenvalue $\lambda$  of
$(D^{1/3})^2$ satisfies the following estimate depending on
the ratio $c:=\Scalg_{\min} / ||T||^2>0$:
\bdm
\lambda \ \geq\ \left\{\ba{ll} \frac{1}{4}\left[\Scalg_{\min}-
\frac{1}{2}||T||^2\right] \, = \, 
\frac{||T||^2}{4}\left[c -\frac{1}{2}\right]
& \text{for } c\geq 3/2 \\[2mm] 
\frac{1}{16}\left[\sqrt{6\,\Scal^g_{\min}} - ||T||\right]^2\, = \,
\frac{||T||^2}{16}\left[\sqrt{6\,c}-1\right]^2 & \text{for } 
1/6\leq c\leq 3/2   \ea\right.
\edm
For $c=3/2$, both estimates coincide. For $c<1/6$, no lower bound can 
be given.
\end{thm}
\begin{NB}
%---------
Hence, for $c\geq 3/2$, the universal bound is still the best one, while
we can improve it for $c\in [1/2,3/2]$ and obtain a new estimate
in the range $c\in [1/6,1/2]$, where the universal bound is negative.
A graph of these estimates is given in Figure \ref{new-estimate} (solid line).
The dashed line shows the values of the other estimate in the interval
where it is not applicable. In particular, we obtain no
 estimate for $c\leq 1/6$, though the curve corresponding to the
$S$-deformed estimate would be positive. Recall that $c=3/2$ 
corresponds to vanishing $\nabla$-scalar curvature, and that examples
with $\nabla$-parallel spinors and hence $\lambda=||T||^2/4$ do exist.
\end{NB}
%-------------------------------------------------------------------------
%--------------------------------------------------------------------------
%
\begin{proof}
%------------
Let $\psi\in\Sigma_{\pm}$ be an eigenspinor of $(D^{1/3})^2$ with
eigenvalue $\lambda$. Our starting point is the $S$-deformed S-L-formula c) 
from Proposition \ref{general-formulas}. First, an easy computation shows that
\bdm
\langle(D^{1/3}+S)^2\psi,\psi\rangle\ =\ \lambda\,||\psi||^2+
2(a_0\pm a_1||T||)\langle D^{1/3}\psi,\psi\rangle + ||S\psi||^2.
\edm
For the right hand side, choose an orthonormal basis such that
$T$ is proportional to $e_{123}$. Then $e_4\cdot T+ T\cdot e_4=0$, while
$e_i\cdot T=T\cdot e_i$ for $i=1,2,3$. Thus, we can compute
\bea[*]
\lefteqn{\sum_{i=1}^4 || (e_i \cdot S + S\cdot e_i)\psi||^2\ =\
\sum_{i=1}^4||2a_0\,e_i\cdot \psi+ a_1(e_i\cdot T+T\cdot e_i)\cdot\psi||^2
}\\
 & =&  ||2\, a_0\,e_4\cdot\psi||^2 + \sum_{i=1}^3 ||2\, a_0\,e_i\cdot\psi +
2\, a_1\, e_i\cdot T\cdot \psi||^2\\
&=&
4\,a_0^2 ||e_4\cdot\psi||^2 + 4 \sum_{i=1}^3 || (a_0\pm a_1 ||T||)e_i
\cdot\psi||^2\ =\ 4 \left[ 3\,(a_0\pm a_1||T||)^2 +\,a_0^2\right]\,||\psi||^2.
\eea[*]
Thus, evaluation of the $S$-deformed S-L-formula yields
\begin{eqnarray*}
\lambda ||\psi||^2 + 2 (a_0\pm a_1||T||)
\langle D^{1/3}\psi,\psi\rangle &=&
 ||\nabla^S\psi||^2 - \left[ 3\,(a_0\pm a_1||T||)^2 +\,a_0^2\right]\,
||\psi||^2 \\
&+& \left[\frac{1}{4}\Scalg -\frac{1}{8}||T||^2 - (a_0\pm a_1||T||)||T||
 \right]\,||\psi||^2.
\end{eqnarray*}
We divide by $||\psi||^2$, exploit that $||\nabla^S\psi||^2/||\psi||^2 
\geq 0$ and
introduce the abbreviation $y:=\langle D^{1/3}\psi,\psi\rangle/||\psi||^2$,
\bdm
\lambda \ \geq\ - (a_0\pm a_1||T||)(2y +||T||)-4\,a_0^2-3 a_1^2||T||^2
\mp 6\, a_0a_1||T|| + \frac{\Scalg_{\min}}{4} 
-\frac{1}{8}||T||^2 .
\edm
This is the rough $S$-deformed eigenvalue estimate which we now need
to concretize. We search the 
global maximum with respect to the variables $a_0,a_1$. 
As partial derivatives, we obtain
\bdm
-(2y +||T||)-8\,a_0\mp 6\,a_1||T||\ = \ 0 \, , \quad 
\mp(2y +||T||)||T|| - 6\,a_1||T||^2\mp 6\, a_0||T|| \ = \ 0 \, .
\edm
The vanishing of both expressions is only possible for
 $a_0=0$ or $T=0$. If $T=0$, the solution $a_0=-y/4$ yields Friedrich's
inequality, see \cite{Friedrich80}. For the case we are interested in, 
$T\neq 0$, we conclude 
$a_0=0$ and  $a_1=\mp (2y+||T||)/6||T||$. Thus, the eigenvalue estimate
becomes
\bdm
\lambda \ \geq\ \frac{1}{12}(2y +||T||)^2 
+ \frac{1}{4}\Scalg_{\min} -
\frac{1}{8}||T||^2.
\edm
At this point, the distinction between $\Sigma_+ $ and $\Sigma_-$
disappears. The Cauchy-Schwarz inequality
implies $y^2\leq \lambda$, i.\,e.~we are interested in the minimum
of $(2y +||T||)^2$ on the interval $[-\sqrt{\lambda},\sqrt{\lambda}]$.

\noindent
\textbf{Case 1:} $\lambda\geq ||T||^2/4$. In this situation, 
the minimum of $(2y +||T||)^2$ on $[-\sqrt{\lambda},\sqrt{\lambda}]$
vanishes, hence the estimate is reduced to the universal bound known before.

\noindent
\textbf{Case 2:} $\lambda <||T||^2/4$. The minimum of $(2y +||T||)^2$
is realized at $y=-\sqrt{\lambda}$, where it takes the value
$(-2\sqrt{\lambda}+||T||)^2= 4\lambda -4\sqrt{\lambda}||T|| +||T||^2$.
Thus we obtain
\bdm
\lambda \ \geq\ \frac{\lambda}{3} -\frac{\sqrt{\lambda}||T||}{3}+
\frac{1}{4}\Scalg - \frac{1}{24}||T||^2,\quad \text{i.\,e. }\,
2\lambda +\sqrt{\lambda}||T|| + \frac{||T||^2}{8}-\frac{3}{4}\Scalg\ \geq\ 0.
\edm
This quadratic inequality yields that  $\sqrt{\lambda}\geq (\sqrt{6\,\Scalg}-||T||)/4$.
For $c\geq 1/6$, the right hand side is positive and we obtain
\bdm
\lambda\ \geq\ \frac{1}{16}(\sqrt{6\,\Scalg}-||T||)^2.
\edm
Now let us discuss which estimate is to be taken for different values of $c$.
We have $\lambda\geq ||T||^2/4$ in the first case, and
$\lambda\geq ||T||^2(\sqrt{6\,c}-1)^2/16$ in the second case;
but there is no way of knowing which case is realized, hence we have to
take the minimum of both values. 

\noindent
If $c\in [1/6,3/2]$, the inequality $1/4 \geq (\sqrt{6\,c}-1)^2/16$ implies 
that we can only conclude $\lambda\geq ||T||^2(\sqrt{6\,c}-1)^2/16$. Since, 
on the other side, this estimate is better than the universal one on 
$[1/6,3/2]$, the claim for this interval follows.

\noindent
For $c\geq 3/2$, $1/4 \leq (\sqrt{6\,c}-1)^2/16$, so we were to conclude
that $\lambda\geq ||T||^2/4$, obviously a rather bad estimate. In this
situation, the universal bound is better and should hence be taken.
\end{proof}
%
%\begin{NB}
%%---------
%A careful inspection of the proof shows that the estimates
%from Theorem \ref{thm-n4}
%still hold for a compact $4$-dimensional spin manifold $(M^4,g)$
%with $\Scalg_{\min}>0$ and a \emph{closed $3$-form $0\neq T\in \Lambda^3(M)$
%of constant length}. This is of course a larger class of manifolds than
%those with parallel $T$.
%\end{NB}
%
\begin{NB}[limiting case]
%-----------------------
If the lower bound of Theorem \ref{thm-n4} ($T \neq 0$)
is an eigenvalue and $\psi$ the eigenspinor, there exists an endomorphism 
$S = a_1 T$ such that $\nabla^S \psi = 0$, i.\,e.
\bdm
\nabla^g_X\psi\ =\ - (a_1+1/4) (X\haken T)\cdot \psi.
\edm
Observe that from the proof, $a_1$ is not an arbitrary constant, but depends 
implicitly on $\psi$ and the action of $D^{1/3}$ on it. The
$\nabla^g$-parallel vector field
$* T$ splits the universal covering $\tilde{M} = N^3 \times \R^1$ and
$\psi$ is a real Killing spinor on the $3$-dimensional manifold 
$N^3$, $\nabla^{N^3}_X = a_1 \, ||T|| \, X \cdot
\psi$. Consequently, $N^3$ is either flat ($a_1 = 0$) 
or isometric to a sphere ($a_1 \neq 0$) (see \cite{Friedrich80}, 
\cite{Dirac-Buch}, \cite{Alexandrov&G&I01}).
\end{NB}
\begin{NB}
%---------
By using the generalized Casimir operator, it was shown in
 Proposition~3.5 of   \cite{Agricola&F04b} that 
$\Scalg_{\min}\geq \frac{3}{16}||T||^2$, i.\,e.~$c\geq 3/16$ implies
that the operator $D^{1/3}$ has  trivial kernel. Theorem \ref{thm-n4}
shows that this conclusion holds even for all $c > 1/6$. Notice
that the Casimir operator can also be used to derive an eigenvalue
estimate for $(D^{1/3})^2$, but it typically stays below the universal bound,
hence is rather useless. For example, in the $4$-dimensional
case considered here it yields 
$\lambda\geq \frac{1}{8}(\Scalg_{\min}+||T||^2/2)$ for $c\geq 3/2$
\cite[Prop. 3.2, Prop. 3.3]{Agricola&F04b}.
\end{NB}
\begin{NB}
%---------
The condition $\nabla T=0$ implies $\nabla^g *T=0$, i.\,e.~there
exists an LC-parallel $1$-form on $(M^4,g)$. For such manifolds,
it has been shown by Alexandrov, Grantcharov and Ivanov in 
\cite{Alexandrov&G&I98} that
the eigenvalues of the Riemannian Dirac operator are bounded by
\bdm
\lambda\big((D^g)^2 \big) \ \geq \ \frac{3}{8} \, \Scalg_{\min}.
\edm
However, $(D^{1/3})^2 = (D^g)^2 + (T D^g + D^g T)/4 + T^2/16$ 
is a perturbation of $(D^g)^2$ by an unbounded operator, i.\,e., the
spectra of $(D^g)^2$ and  $(D^{1/3})^2$ are not related in any obvious manner.

In \cite{Alexandrov&G&I01}, the same authors proved that for 
the Dolbeault operator on a Hermitian spin surface with strictly positive
conformal scalar curvature $k$, the estimate
\bdm\tag{$*$}
\lambda^2 \ \geq\ \frac{1}{2}k_{\min}
\edm
holds. Under our assumption $\nabla T=0$, the Lee form $\theta$ is coclosed,
hence the relation between Riemannian scalar and $*$-scalar curvature reads
$\Scalg-\Scal^*=||\theta||^2=||T||^2$. By definition, $k=(3\Scal^*-\Scalg)/2$,
so $k>0$ is equivalent to $c>3/2$ and the estimate ($*$) is
equivalent to $\lambda^2\geq ||T||^2(2c-3)/4$. This is a line
going through zero at $c=3/2$ and that becomes better than the universal
estimate for $c\geq 5/2$. We think that it cannot be derived in the
more general framework described here (where no Hermitian structure is
assumed).
\end{NB}
\begin{exa}
%----------
Consider the $2$-dimensional sphere $X^2$ equipped with a
Riemannian metric of positive Gaussian curvature $G$. 
Denote by
$G_{\min}$ its minimum.
Moreover, we fix a positive number $||T||$ such that the following conditions
hold:
\bdm
|| T ||^2 \ < \ G_{\min} \, , \quad \frac{||T|| \, \mathrm{vol}(X^2)}{\pi} \ =
\ k \ \ \text{is an integer} \, .
\edm
The $2$-form $F := 2 ||T|| \, dX^2$ satisfies the condition
\bdm
\frac{1}{ 2 \pi} \int_{X^2} F \ = \ k  .
\edm
Consequently, there exists an $S^1$-principal bundle $N^3 \rightarrow X^2$ and
a connection form $\eta$ on $N^3$ such that $d \eta = \pi^*(F)$. We split
the tangent bundle $TN^3 = T^v \oplus T^h$ into its vertical and horizontal
part and we define a metric $g$ on $N^3$ by pulling back the metric of the
surface $X^2$. The complex structure of $X^2$ lifts to an endomorphism
$\varphi : TN^3 \rightarrow TN^3$ such that
\bdm
\varphi^2 \ = \ - \, \mathrm{Id} \, + \, \eta \otimes \eta \, , \quad
g(\varphi(X) , \varphi(Y)) \ = \ g(X,Y) \, - \, \eta(X) \eta(Y) \, , \quad
d \eta \ = \ 2 ||T|| \, \pi^*(dX^2)
\edm
holds. The tuple $(N^3, g, \eta, \varphi)$ is a $3$-dimensional compact
$\alpha$-Sasakian manifold with fundamental form $\pi^*(dX^2)$.
The scalar curvature
$\Scalg$ of $N^3$ is given by the formula
\bdm
\Scalg \ = \ 2 \, G \ - \ 2 \, ||T||^2 \ > \ 0 \, .
\edm
The Riemannian product $M^4 := N^3 \times S^1$ admits a canonical
complex
structure $J$ such that the K\"ahler form is given by the formula
$\Omega =  dt \wedge \eta +  \pi^*(dX^2)$.
Consider the $1$-form
$\Gamma := ||T|| dt$ . Then we obtain
\bdm
d \Omega \ = \ - \, dt \wedge d \eta \ = \ - \, 
2 \, ||T|| \, dt \wedge \pi^*(dX^2)
\ = \ - \, 2 \, \Gamma \wedge \Omega \, .
\edm
Consequently, $(M^4, g \oplus dt^2,J)$ is a compact
Hermitian $4$-manifold with positive scalar curvature and parallel
characteristic torsion (see \cite{Agricola&F06} and \cite{Alexandrov&F&S04}),
\bdm
T^c \ = \ - \, ||T|| \, dN^3 \, , \quad
\nabla^c \ = \ \nabla^g \, + \, \frac{1}{2} T^c \, , \quad \nabla^c \Omega \ =
\ 0 \, , \quad \nabla^c T^c \ = \ 0 \, .
\edm
By the way, any Hermitian $4$-manifold with parallel characteristic torsion
and positive scalar curvature is locally isometric to a manifold
of our family. Now,
the universal lower bound $\beta_{\mathrm{univ}}$, the $S$-deformed lower 
bound $\beta_{S}$ for $(D^{1/3})^2$ on $M^4$
as well as the relevant ratio $c$
are given by
\bdm
\beta_{\mathrm{univ}}\ =\ \frac{1}{2} G_{\min}  -  \frac{5}{8}||T||^2, \quad
\beta_{S}\ =\ \frac{1}{16}\big(\sqrt{12(G_{\min}-||T||^2)}- ||T||\big)^2,\quad 
c \ = \ \frac{2 G_{\min}}{||T||^2}  -  2  .
\edm
%
%-------------
\begin{figure}
\noindent
%--------------------------------
\begin{minipage}[b]{.46\linewidth}
\centering
\psfrag{2}{$2$}\psfrag{a}{$a$}\psfrag{c}{$c$}
\psfrag{10}{$10$}\psfrag{8}{$8$}\psfrag{6}{$6$}\psfrag{4}{$4$}
\psfrag{1/4}{$1/4$}\psfrag{1/2}{$1/2$}\psfrag{3/4}{$3/4$}
\psfrag{1}{$1$}\psfrag{3/2}{$3/2$}\psfrag{5/4}{$5/4$}
\includegraphics[width=\linewidth]{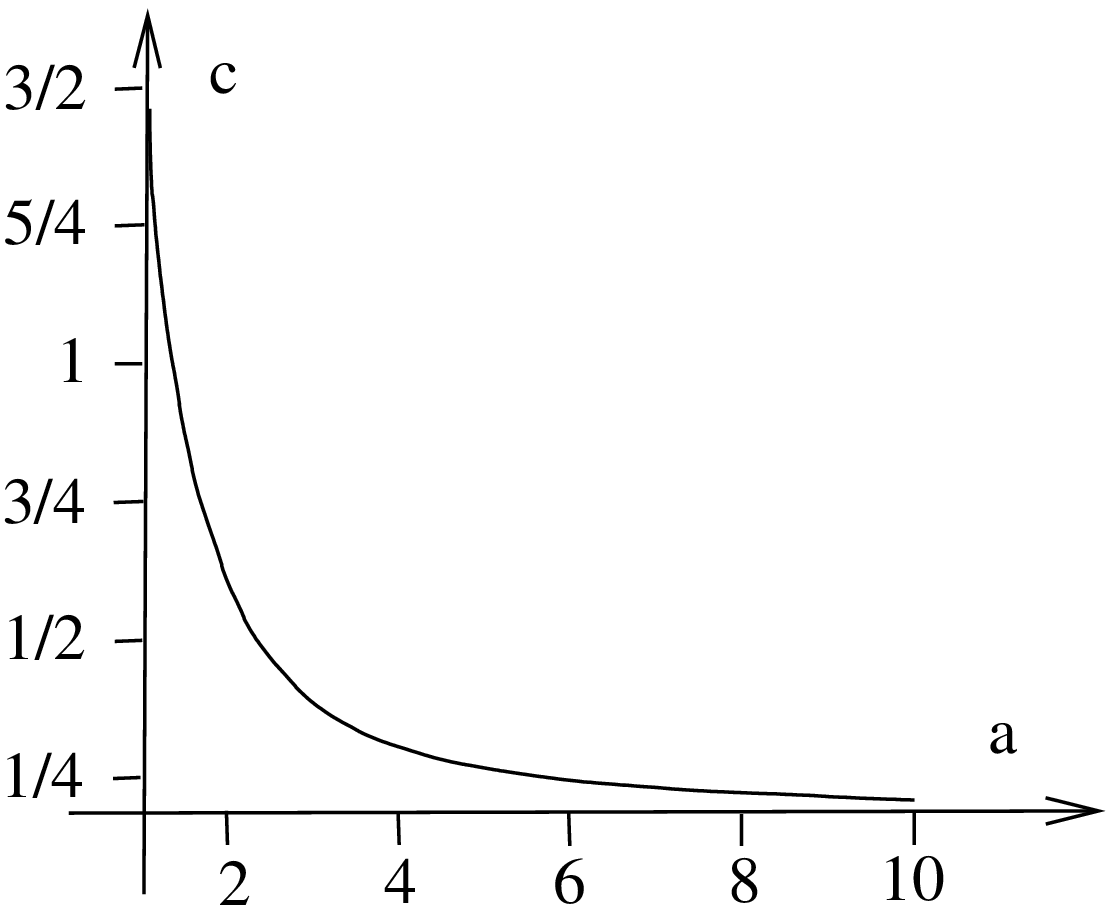}
\caption{The parameter $c=\Scalg_{\min}/||T||^2$ for $a>1$.}
\label{c-for-a-large}
\end{minipage}\hfill
%---------------------------------
\begin{minipage}[b]{.48\linewidth}
\centering
\psfrag{a}{$a$}
\psfrag{0.01}{$0.01$}\psfrag{0.02}{$0.02$}\psfrag{0.03}{$0.03$}
\psfrag{0.005}{$0.005$}\psfrag{0.015}{$0.015$}\psfrag{0.025}{$0.025$}
\psfrag{1}{$1$}\psfrag{3/2}{$3/2$}\psfrag{2}{$2$}
\psfrag{3}{$3$}\psfrag{4}{$4$}\psfrag{5}{$5$}
\psfrag{Bs}{$\beta_S$}\psfrag{Bu}{$\beta_{\mathrm{univ}}$}
\includegraphics[width=\linewidth]{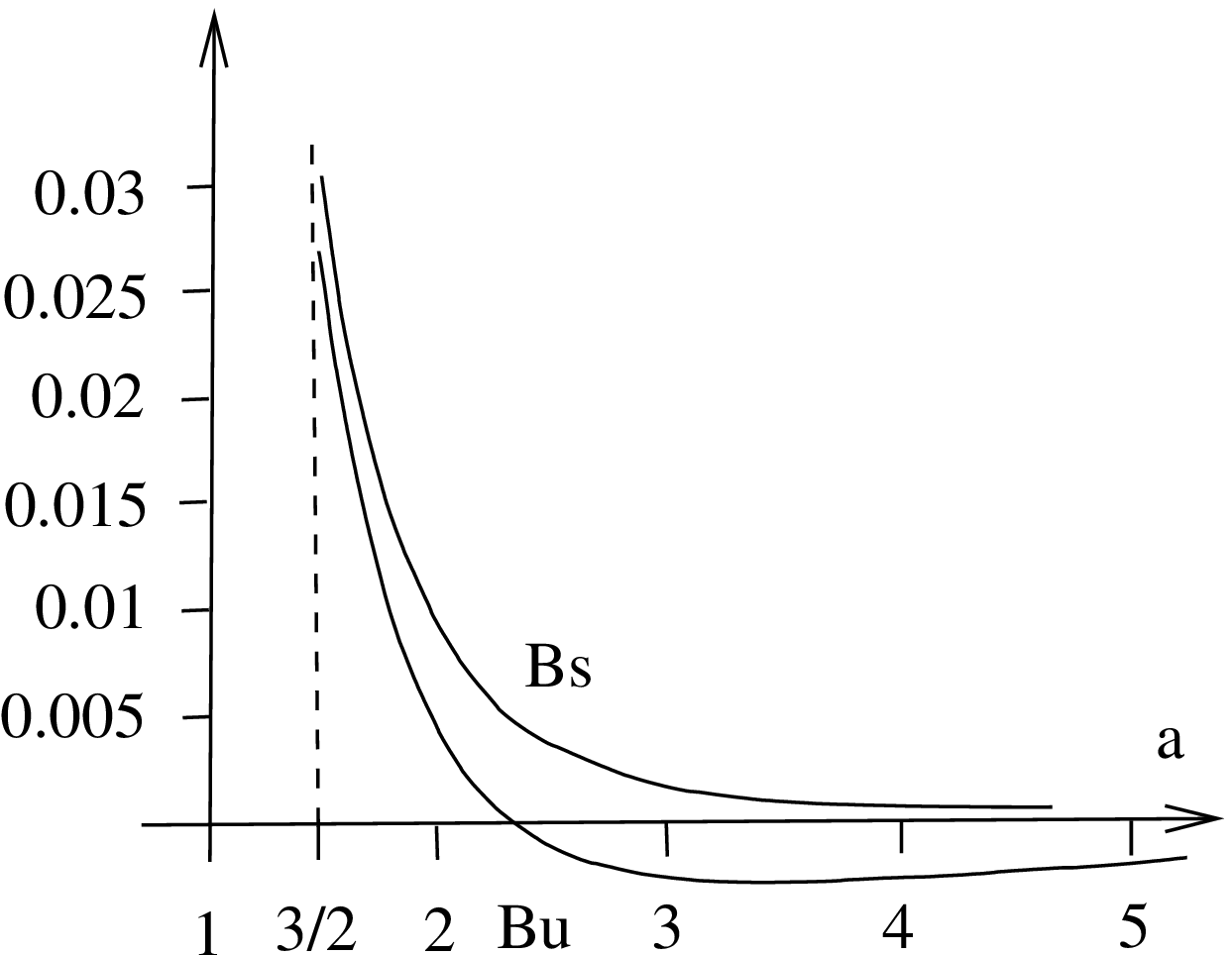}
\caption{Lower bounds $\beta_{\mathrm{univ}}$ and $\beta_{S}$ for $a>1$.}
\label{bounds-for-a-large}
\end{minipage}
\end{figure}
%--------------------------------------------------------------------
%
Let us discuss these estimates in case of an ellipsoid with the induced metric
(see for example \cite{Agricola&F99} for the relevant formulas)
\bdm
X^2 \ = \ \big\{ (x,y,z) \in \R^3 \, : \, x^2  + y^2  +
\frac{z^2}{a^2} \ = \ 1 \big\} \, .
\edm
Computing the minimum of the Gaussian curvature and the volume, 
\bdm
G_{\min}=\frac{1}{a^2},\quad \vol(X^2)\ =\ 
2\pi+ 2\pi a^2 \, \frac{\arcsin\left[\frac{1}{a}\sqrt{a^2-1}\right]}{
\sqrt{a^2-1}},
\edm
we see that the
conditions $|| T ||^2 <  G_{\min}$ and  $||T|| \, \mathrm{vol}(X^2) =
k \, \pi$ imply that there are only three admissible values for $k$, namely
$k = 1,2,3$.
If $a >1$,
the ratio $c=\Scalg_{\min}/||T||^2$ is always larger than $1/6$ 
and $c < 3/2$ is true for
$1.02 < a < \infty$ (see Figure \ref{c-for-a-large} for $k=3$), although 
$\lim_{a\ra 1^+} c=14/9= 1.555\ldots$. 
Figure \ref{bounds-for-a-large} shows the universal bound 
$\beta_{\mathrm{univ}}$ (lower curve) and the new bound $\beta_S$
(upper curve) as a function of $a$ for $a\geq 3/2$. While 
$\beta_{\mathrm{univ}}(a)$  becomes negative for $a>2.33$, the new bound 
$\beta_S(a)$ stays strictly positive for all
$a$ (and has limit $0$ when $a\ra\infty$). For $a\ra 1^+$, the two curves 
approach each other more and more, though 
$\beta_S(a)>\beta_{\mathrm{univ}}(a)$ holds for all $a$. Their limits
are
\bdm
\lim_{a\ra 1^+} \beta_{\mathrm{univ}} \ =\ \frac{19}{128}\ 
\cong\ 0.1484375,\quad
\lim_{a\ra 1^+} \beta_S \ =\ \frac{93}{256} - \frac{3\sqrt{21}}{64}\ \cong\
0.1484730143.
\edm
Consequently, the $S$-deformed lower bound on $M^4$ is better than the 
universal bound for ellipsoidal deformations of the sphere  with $a>1.02$. 
A similar discussion applies in case of $a < 1$.

\end{exa}
\begin{NB} The example also shows  that the $S$-deformed lower bound 
$\beta_S$ applies \emph{only}
under the condition that $c \leq 3/2$. Indeed, consider the (forbidden) case 
$||T|| = 0$. Then $M^4 = X^2 \times T^2$ is the Riemannian product
of $X^2$ by the flat torus $T^2$. The scalar curvature $\Scalg$ of $M^4$
equals
$2G$ and $D^{1/3}$ coincides with the Riemannian Dirac operator $D^g$ of
$M^4$. The $S$-deformed lower bound for $(D^{1/3})^2$ yields $3 \Scalg/2$. 
However, the sharp estimate of the Dirac operator on a $4$-dimensional
K\"ahler manifold is $\Scalg/2$ and this lower bound is realized 
on $S^2 \times T^2$ (see \cite{Kirchberg}, \cite{Friedrich93}).
\end{NB}
%

%
%---------------------------------------------------------------------------- 
\section{$5$-dimensional Sasakian manifolds}
\label{fam-conn}\noindent
%----------------------------------------------------------------------------
%
Let $(M^5,g,\xi,\eta,\varphi)$ be a compact 5-dimensional Sasakian manifold
with a fixed spin structure. There exists a unique connection $\nabla$
with totally
skew-symmetric torsion and preserving the Sasakian structure (see 
\cite{Friedrich&I1}). The torsion form is given by the formula
$T = \eta \wedge d \eta$, $||T||=8$. $T$ splits the spinor bundle into two
$1$-dimensional bundles and one $2$-dimensional bundle,
\bdm
\Sigma_{\pm4} \ = \ \big\{ \psi \in \Sigma M^5 \, : \, T\psi \, = \, \pm 
\, 4 \, \psi \big\} \, , \quad
\Sigma_{0} \ = \ \big\{ \psi \in \Sigma M^5 \, : \, T\psi \, = \, 0 \big\} 
\, .
\edm
In particular, we obtain $T^3 = 16 \, T$ and 
any admissible polynomial in $T$ is
quadratic, $S = P(T) = a_0 I \, + \, a_1 T \, + \, a_2 T^2$. We will estimate
$(D^{1/3})^2$ in $\Sigma_{\pm4}$ and in  $\Sigma_{0}$ separately.
%
%---------------------------------------------------------------------------- 
\subsection{The estimate of $(D^{1/3})^2$ in the subbundle $\Sigma_{0}$}
\label{fam-conn}\noindent
%----------------------------------------------------------------------------
%
\noindent
The operator  $D^{1/3} =D^g + T/4$ coincides in the subbundle
$\Sigma_0$ with
the Riemannian Dirac operator. Consequently, we can use the estimate
of the Riemannian Dirac operator (see \cite{Friedrich80}),
\bdm
\lambda_{\min} \big((D^{1/3})^2_{|\Sigma_{0}}\big)\ 
\geq \ \lambda_{\min}\big((D^{g})^2\big) \ \geq \ \frac{5}{16}
\Scalg_{\min}\ =:\ \beta^g .
\edm
If the scalar curvature is small, there is a better estimate. 
Indeed, the universal S-L-formula for $(D^{1/3})^2$ on spinors
stated in equation (\ref{D/3-square}) reduces on $\Gamma(\Sigma_{ 0})$ to
\bdm
 (D^{1/3})^2 \ = \ \Delta \, + \, 1 + \frac{1}{4} \Scalg \, .
\edm
Both inequalities together yield the following estimate
\bdm
\lambda_{\min}\big((D^{1/3})^2_{|\Sigma_{0}}\big)\ 
\geq \ \max \left\{ \frac{5}{16} \Scalg_{\min}  , \, 
1  +  \frac{1}{4} \Scalg_{\min} \right\} .
\edm
In case of  $-4 = \Scalg_{\min}$, any spinor in the kernel of 
$(D^{1/3})^2$ is parallel and the Sasakian manifold is isometric
to a compact quotient of the $5$-dimensional Heisenberg group, see
\cite{Friedrich&I2}. If  $-4 < \Scalg_{\min}$, 
the eigenvalues of $(D^{1/3})^2$ on spinors in  $\Gamma(\Sigma_{ 0})$ are
positive and bounded by the eigenvalues on spinors in  
$\Gamma(\Sigma_{\pm 4})$. More precisely, we prove the following
\begin{prop}\label{Sigma0-Sigma4}
%--------------------------------
Let $(M^5,g,\xi,\eta,\varphi)$ be a compact $5$-dimensional
Sasakian manifold such that $-4<\Scalg_{\min}$ holds. If  
$\psi\in\Gamma(\Sigma_0)$ is an eigenspinor of the operator 
$(D^{1/3})^2$ and $\lambda$ is the eigenvalue, then any part of
the decomposition
\bdm
D^{1/3} \psi \ = \ \alpha_4 \, + \, \alpha_0 \, + \, \alpha_{-4} \, , \quad
\alpha_{k} \in \Gamma(\Sigma_{k}) \, , \quad
k\in\left\{-4,0,4\right\} 
\edm
is an eigenspinor of the operator 
$(D^{1/3})^2$ with the same eigenvalue $\lambda$. 
Moreover, at least one of the spinors  $\alpha_{\pm 4}$ is nontrivial.
In particular,  the eigenvalues of $(D^{1/3})^2$ on $\Gamma(\Sigma_{ 0})$ 
are bounded by the eigenvalues of 
$(D^{1/3})^2$ on $\Gamma(\Sigma_{\pm 4})$,
\bdm
\lambda_{\min}\big((D^{1/3})^2_{|\Sigma_{0}}\big)\ \geq\ 
\lambda_{\min}\big((D^{1/3})^2_{|\Sigma_{\pm 4}}\big).
\edm
\end{prop}
\begin{proof}
%------------
The operator $(D^{1/3})^2$ preserves the splitting of the spinor bundle.
Therefore, the components $\alpha_k$ are eigenspinors of  $(D^{1/3})^2$.
If $\alpha_{\pm4}$ vanish for any eigenspinor in $\Sigma_0$, the 
operator $D^{1/3}$ acts on the corresponding eigenspace as a symmetric
operator. Consequently, there exists a spinor field $\psi$ such that
\bdm
D^{1/3} \psi \ = \ \pm \, \sqrt{\lambda} \, \psi \, , \quad 
\psi \, \in \Gamma(\Sigma_0)  .
\edm
By Proposition \ref{polynom}, the connection $\nabla^S$ defined by the
polynomial $S = a ( - 8\,  \mathrm{Id}+ T^2)$ preserves the bundle
$\Sigma_0$. We compute that
\bdm
\sum_{i=1}^5||(e_iS \, + \, Se_i)\psi||^2 \ = \ 2\cdot16^2a^2 \, ||\psi||^2
\edm
holds for any spinor $\psi \in \Sigma_0$. The $S$-deformed S-L-formula
of Proposition \ref{general-formulas} yields the following inequality
\bdm
\lambda \ \geq \
16a\frac{\langle\psi,D^{1/3}\psi\rangle_{L^2}}{||\psi||^2}\, - \, 
\frac{1}{4}16^2a^2 \, + \, 1 \, + \, \frac{1}{4}\Scalg_{\textnormal{min}} \, .
\edm
The optimal parameter 
\bdm
a \ = \ \frac{1}{8}\frac{\langle\psi,D^{1/3}\psi\rangle_{L^2}}{||\psi||^2}
\edm
implies the inequality
\bdm
\lambda \ \geq \ y^2 \, + \, 1 \, + \, \frac{1}{4}\Scalg_{\textnormal{min}} \, ,
\edm
where  $y:=\langle\psi,D^{1/3}\psi\rangle_{L^2}/||\psi||^2 = 
\pm \sqrt{\lambda}$. Finally
we conclude that $-4 \geq \Scalg_{\mathrm{min}}$, a contradiction.
\end{proof}
%
%---------------------------------------------------------------------------- 
\subsection{The estimate of $(D^{1/3})^2$ in the subbundle $\Sigma_{\pm4}$}
\noindent
%----------------------------------------------------------------------------
%
Let us consider the subbundle
$\Sigma_4=\left\{\psi\in\Sigma;\ T\psi=4\psi\right\}$. Then
$TM^5\cdot\Sigma_{4}$ is a proper subbundle of $\Sigma$. Indeed, we have
\bdm
        TM^5\cdot \Sigma_4\ \subset \  \Sigma_0\op \Sigma_4  .
\edm
By Proposition \ref{polynom}, the connection $\nabla^S$ defined by the
polynomial $S = (-2a_1-8a_2)\Id +a_1 T+a_2 T^2$
 preserves the bundle $\Sigma_4$. In particular, $S$ acts 
in $\Sigma_4$ by multiplication with $x = 2a_1+8a_2$. 
The operator  $(D^{1/3})^2$ on spinors in  $\Gamma(\Sigma_{ 4})$ is given
by the formula (see eq.~(\ref{D/3-square}))
\bdm
(D^{1/3})^2_{|\Sigma_4}\ = \ \Delta \, -3 \, + \, \frac{1}{4}\Scalg 
\edm
and we obtain the following universal estimate for  all eigenvalues $\lambda$
of the operator $(D^{1/3})^2_{|\Sigma_4}$,  
\bdm
\lambda\ \geq \ -3+\frac{1}{4}
\Scalg_{\min}\ =:\ \beta_{\mathrm{univ}}  .
\edm
Proposition \ref{Sigma0-Sigma4} implies that this inequality
also holds for all eigenvalues of $(D^{1/3})^2$. In order to
evaluate the $S$-deformed eigenvalue bound, we first 
compute  for spinors $\psi\in \Sigma_4$
\bdm
\sum_{i=1}^5||(e_iS+Se_i)\psi||^2 \ = \, 16(a_1+4a_2)^2 ||\psi||^2 .
\edm
The $S$-deformed S-L-formula
of Proposition \ref{general-formulas} then yields the  inequality
\bdm
\lambda \ \geq \ 
-x^2-4x-2x\frac{\langle\psi,D^{1/3}\psi\rangle}{||\psi||^2}-
3+\frac{1}{4}\Scalg_{\textnormal{min}} \, , \quad x\ :=\ 2a_1+8a_2  .
\edm
The optimal parameter
\bdm
x \ = \ -2-\frac{\lan\psi,D^{1/3}\psi\ran_{L^2}}{||\psi||_{L^2}^2}
\edm
yields the refined estimate
\bdm
\lambda \ \geq \ 1+4y+y^2+\frac{1}{4}
\Scalg_{\textnormal{min}} \, , \quad
y\ :=\ \lan\psi,D^{1/3}\psi\ran_{L^2}/||\psi||^2 \, .
\edm
Let us discuss the limiting case of $\lambda = 1+ 4 y + y^2 + 
\frac{1}{4}\, \Scalg_{\mathrm{min}}$. Then the scalar curvature
is constant and the eigenspinor $\psi \in \Gamma(\Sigma_{\pm 4})$
is parallel with respect to a connection given by the formula
\bdm
\nabla^S \psi \ = \ 0, \quad S \ = \ (- 2a_1 - 8 a_2)\Id+ a_1T+a_2 T^2.
\edm
Conversely, an easy algebraic computation yields that any $\nabla^S$-parallel
spinor in $\Sigma_{\pm 4}$ is an eigenspinor of the operator $D^{1/3}$ with
eigenvalue $y = - (2 + 2 a_1 \pm 8 a_2)$. In the limiting case
$y^2 = \lambda = 1 + 4 y + y^2 + \frac{1}{4}\, \Scalg_{\mathrm{min}}$ we
obtain
\bdm
\lambda \ = \ y^2 \ = \ \frac{1}{16} \left[ 1 \, + \, \frac{1}{4}
\Scalg_{\min} \right]^2  .
\edm
\begin{prop}
%-----------
For any eigenspinor $\psi \in \Gamma(\Sigma_{\pm 4})$ with eigenvalue
$\lambda$ holds
\bdm
\lambda \ \geq \ 1+4y+y^2+\frac{1}{4}
\Scalg_{\textnormal{min}} \, , \quad
y:=\lan\psi,D^{1/3}\psi\ran_{L^2}/||\psi||^2 \, .
\edm 
Equality occurs if and only if the scalar curvature is constant and there
exist numbers $a_1, a_2$ such that
\bdm
\nabla^S \psi \ = \ 0 \, , \quad S \ = \ (- 2a_1 \, - \, 8 a_2) I \, + \, a_1
T \, + \, a_2 T^2 \, , \quad \Scalg \ = \ 28 \, + \, 32 a_1 \, \pm \, 
128 a_2  .
\edm
In this case we have
\bdm
\lambda \ = \ y^2 \ = \ \frac{1}{16} \left[ 1 \, + \, \frac{1}{4}
\Scalg_{\min} \right]^2 \, .
\edm
\end{prop}
\noindent
Since $\Sigma_{\pm 4}$ are $1$-dimensional bundles, the integrability
condition
for the existence of a parallel spinor, $\nabla^S \psi = 0$, can be 
formulated using the Ricci tensor only. A standard argument 
(see \cite{FriedrichKim}) yields
the following characterization.
\begin{prop}
%-----------
A $5$-dimensional simply connected Sasakian manifold admits a spinor 
field $\psi \in\Gamma(\Sigma_{\pm4})$ such that
\bdm
\nabla^S \psi \ = \ 0 \, , \quad S \ = \ (- 2a_1 \, - \, 8 a_2) I \, + \, a_1
T \, + \, a_2 T^2 \, ,
\edm
if and only if it is $\eta$-Einstein. The Riemannian Ricci tensor
has the eigenvalues  $6 + 8 a_1 \pm 32 a_2$ with multiplicity four and
$4$ with multiplicity one.
\end{prop}
\begin{cor}
%-----------
$\frac{1}{16} \big[ 1 +  \frac{1}{4}
\Scalg \big]^2 $ is an eigenvalue of the operator $(D^{1/3})^2$ on 
any $5$-dimensional $\eta$-Einstein-Sasakian manifold.
\end{cor}
\noindent
Now we estimate the eigenvalues of $(D^{1/3})^2$ by purely geometric data. 
The Cauchy-Schwarz inequality $y^2 \leq \lambda$
restricts the parameter $y$. If $\lambda_{\min}< 4$, we obtain 
\begin{eqnarray*}
\underset{y\in\left[-\sqrt{\lambda},\sqrt{\lambda}\right]}{\min}
\left\{1+4y+y^2\right\}\
= \ \left\{1+4y+y^2\right\}_{|y=-\sqrt{\lambda}}\ = \ 
1-4\sqrt{\lambda}+\lambda\, .
\end{eqnarray*}
Finally, we estimated the operator
$(D^{1/3})^2$ on spinors in $\Gamma(\Sigma_4)$,
\bdm
\lambda_{\min} \ \geq \ \min\left\{4,\frac{1}{16}\left[1+\frac{1}{4}
\Scalg_{\min}\right]^2\right\} \, .
\edm
However, for large scalar curvature $\Scalg_{\min}\geq 28$,
the universal estimate lies above $4$,  
$\beta_{\mathrm{univ}}=- 3 + \Scalg_{\min}/4\geq 4$, hence 
\bdm
\lambda\ \geq \ \left\{\ba{ll}
\frac{1}{16} \left[1+ \frac{1}{4}\Scalg_{\min}\right]^2 & 
\text{ for } -4\ \leq \ \Scalg_{\mathrm{min}}\leq 28,\\[2mm]
-3+ \frac{1}{4}\, \Scalg_{\mathrm{min}} & \text{ for } \Scalg_{\min}\geq 28. 
\ea\right. 
\edm
\noindent
We study the case  $\Scalg_{\min} \geq -4$ more carefully.
Let $\psi \in \Gamma(\Sigma_4)$ be an eigenspinor of $(D^{1/3})^2$ realizing
the minimal eigenvalue. The spinor $D^{1/3}\psi$ is thus a section of 
the bundle $\Sigma_0 \oplus \Sigma_4$.
We decompose $D^{1/3}\psi = \alpha_0 + \alpha_4$ into its components. Then
$\alpha_0$ and $\alpha_4$ are eigenspinors of $(D^{1/3})^2$. If $\alpha_0 \neq
0$ we conclude that
\bdm
 \lambda_{\min}\big((D^{1/3})^2_{|\Sigma_{\pm4}}\big) \ \geq \
 \lambda_{\min}\big((D^{1/3})^2_{|\Sigma_{0}}\big) \, .
\edm
If $\alpha_0 = 0$ for any eigenspinor related to the minimal eigenvalue
$\lambda_{\mathrm{min}}$, the operator $D^{1/3}$ is a symmetric
operator acting in the eigenspace $\{ \psi : (D^{1/3})^2 \psi =
  \lambda_{\mathrm{min}} \psi , \psi \in \Gamma(\Sigma_4)\}$. Consequently,
there exists an eigenspinor of $D^{1/3}$ inside the bundle $\Sigma_4$,
\bdm
D^{1/3} \psi \ = \ \lambda_{\mathrm{min}} \, \psi \, , \quad
\psi \ \in \ \Gamma(\Sigma_4) \, .
\edm
Then we obtain $ y = \pm \sqrt{\lambda_{\mathrm{min}}}$ and the inequality
$ \lambda \geq  1+4y+y^2+\frac{1}{4}
\Scalg_{\min}$ yields the estimate
\bdm
\lambda_{\min} \ \geq \ 
\frac{1}{16}\left[1+\frac{1}{4}\Scalg_{\textnormal{min}}\right]^2 
\ =:\ \beta_S.
\edm 
In particular, we proved
\begin{prop}\label{Sigma4-min}
%-----------------------------
\bdm
\lambda_{\min}((D^{1/3})^2_{\Sigma_{\pm4}}) \ \geq \
\min\left\{\frac{1}{16}\left[1+\frac{1}{4}\Scalg_{\textnormal{min}}\right]^2 ,
\, \lambda_{\min}\big((D^{1/3})^2_{|\Sigma_{0}}\big) \right\} .
\edm
\end{prop}
\noindent
Let us summarize the previous discussion.
The inequalities of Proposition \ref{Sigma0-Sigma4}, Proposition 
\ref{Sigma4-min} as well as the inequality 
$\lambda_{\min}\big((D^{1/3})^2_{|\Sigma_{0}}\big)\geq \frac{5}{16}
\Scalg_{\min}=:\beta^g$ together yield the following result.
\begin{thm}\label{est-dim5}
%--------------------------
Let $(M^5,g,\xi,\eta,\varphi)$ be a  compact Sasakian
manifold with $\Scalg_{\min}>-4$, $T=\eta\hut d\eta$ its characteristic 
torsion. The first eigenvalue $\lambda_{\min}$ of $(D^{1/3})^2$ satisfies:
\bdm
\lambda_{\min}\ \geq\ \left\{\ba{ll} 
\frac{1}{16}\left[1+\frac{1}{4}\Scalg_{\min}\right]^2 & \text{ for } 
-4<\Scalg_{\min}\leq 4(9+4\sqrt{5})\\[2mm]
\frac{5}{16}\, \Scalg_{\min} & \text{ for }\ \Scalg_{\min}\geq 4(9+4\sqrt{5})
\simeq 71,78.
\ea\ \right.
\edm
For $\Scalg_{\min}= 4(9+4\sqrt{5})$, both estimates coincide.
Furthermore, the smallest eigenvalues of the operators
$(D^{1/3})^2_{|\Sigma_{\pm 4}}$ and $(D^{1/3})^2_{|\Sigma_0}$ satisfy:
\begin{enumerate}
\item $\lambda_{\min}((D^{1/3})^2_{|\Sigma_0})  \geq 
\lambda_{min}((D^{1/3})^2_{|\Sigma_{\pm 4}})$ for all $\Scalg_{\min}>-4$,
hence $\lambda_{\min}$ can always be realized in $\Sigma_{\pm 4}$;
\item If, for $\Scalg_{\min}> 4(9+4\sqrt{5})$, $\lambda_{\min}$ happens
to lie in the intermediate range 
$\frac{1}{16}\left[1+\frac{1}{4}\Scalg_{\min}\right]^2\geq \lambda_{\min}\geq
\frac{5}{16} \Scalg_{\min}$, then $\lambda_{\min}((D^{1/3})^2_{|\Sigma_0})
= \lambda_{\min}((D^{1/3})^2_{|\Sigma_{\pm4}})$.
\end{enumerate}
Finally,  $\frac{1}{16}\left[1+\frac{1}{4}\Scalg\right]^2$ is an eigenvalue of 
$(D^{1/3})^2$ on any $\eta$-Einstein-Sasakian manifold.
\end{thm} 
%
%-------------------------------------------------------------------------
\begin{figure}
\bdm
\psfrag{-4}{\scriptsize $-4$}
\psfrag{10}{$10$}\psfrag{20}{$20$}\psfrag{30}{$30$}
\psfrag{40}{$40$}\psfrag{60}{$60$}\psfrag{80}{$80$}
\psfrag{c}{\scriptsize $4(9+4\sqrt{5})$}\psfrag{s}{$\Scal^g_{\min}$}
\psfrag{x}{$*$}
\psfrag{Bu}{$\beta_{\mathrm{univ}}$}
\psfrag{Bg}{$\beta^g$}\psfrag{Bs}{$\beta_S$}
\includegraphics[width=6cm]{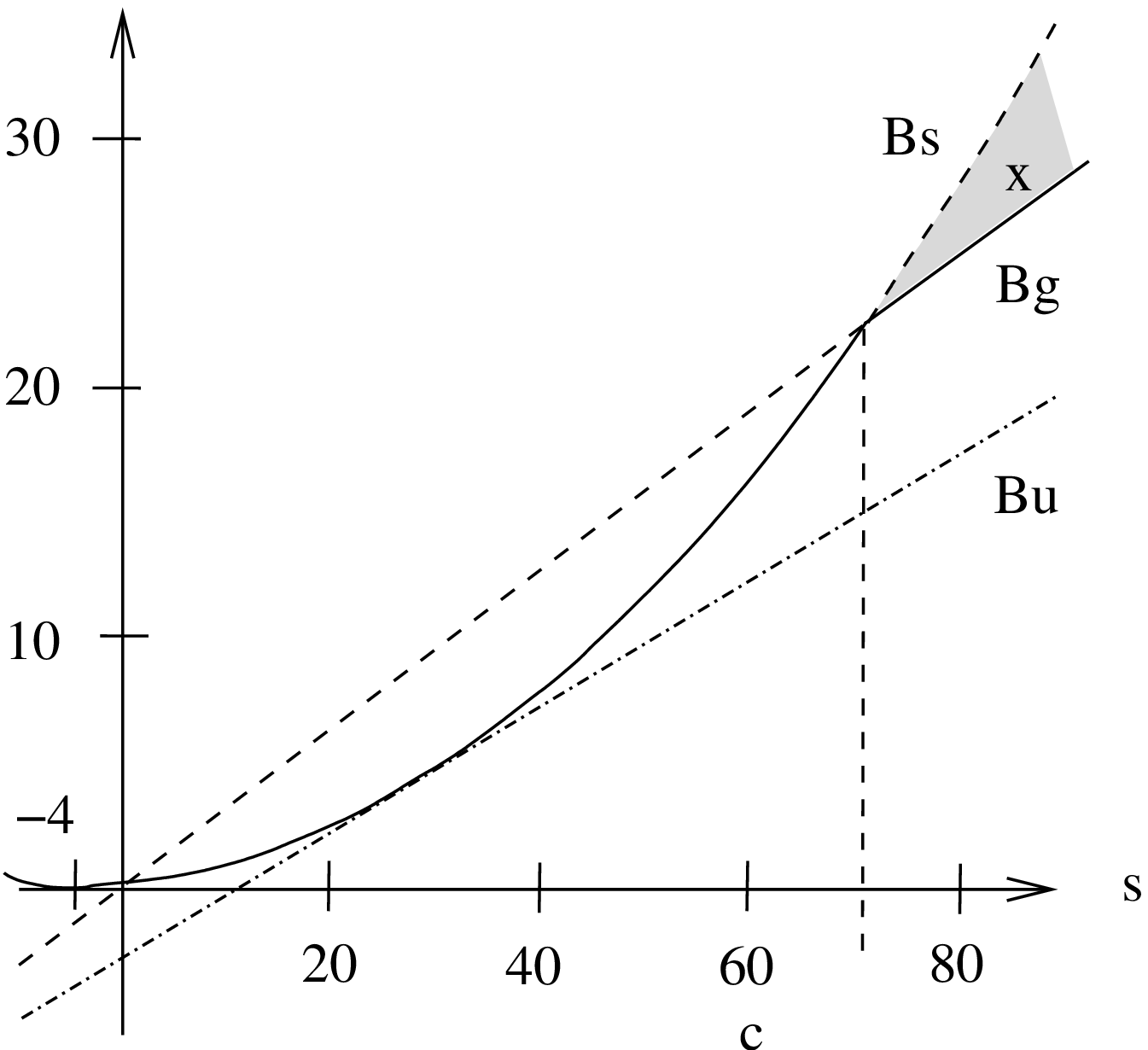}
\edm
\caption{Lower bound for $\lambda$ according to Theorem \ref{est-dim5}.}
\label{new-estimate-5dim}
\end{figure}
%--------------------------------------------------------------------------
\begin{NB}
%---------
A graph of these estimates is given in Figure \ref{new-estimate-5dim}.
The bottom dashed line  corresponds to the universal lower
bound $\beta_{\mathrm{univ}}=-3+\Scalg_{\min}/4$ and turns out to be
useless for all values of $\Scalg_{\min}$. The $S$-deformed lower bound
$\beta_S=\frac{1}{16}\left[1+\frac{1}{4}\Scalg_{\min}\right]^2$ is
drawn as a solid line within its range of application 
$\Scalg_{\min}\in [-4,4(9+4\sqrt{5}]$, and dashed outside.
The Riemannian estimate $\beta^g=5/16\,\Scalg_{\min}$ is to be taken
for $\Scalg_{\min}\geq 4(9+4\sqrt{5})$, corresponding to the solid part of
the line.
A particularly interesting region is the shaded intermediate range $*$
between $\beta^g$ and $\beta^S$; here, $\lambda_{\min}((D^{1/3})^2_{|\Sigma_0})
= \lambda_{\min}((D^{1/3})^2_{|\Sigma_{\pm4}})$ holds.
\end{NB}
%----------------------------------------------------------------------------
%---------------------------------------------------------------------------- 
%------------------------------------------
%\addcontentsline{toc}{section}{Literature}
%------------------------------------------
    

\begin{thebibliography}{1111}
\bibitem{Agricola03}
I. Agricola, \emph{Connections on naturally reductive spaces, their
Dirac operator and homogeneous models in string theory}, Comm. 
Math. Phys. 232 (2003), 535-563.
%
\bibitem{Agricola06}
\bysame, \emph{The Srn\'{\i} lectures on non-integrable geometries with 
torsion}, 
to appear in Suppl. Rend. Circ. Mat. di Palermo Ser. II, 2006.
%
\bibitem{Agricola&F99}
I. Agricola and Th. Friedrich, \emph{Upper bounds for the first
eigenvalue of the Dirac operator on surfaces}, J. Geom. Phys. 30 (1999),
1-22.
%
\bibitem{Agricola&F04a}
\bysame, \emph{On the holonomy of 
connections with skew-symmetric torsion},  Math. Ann.
328 (2004), 711-748.
%
\bibitem{Agricola&F04b}
\bysame, \emph{The Casimir operator of a metric
connection with totally skew-symmetric torsion},
J. Geom. Phys. 50 (2004), 188-204.
%
\bibitem{Agricola&F06}
\bysame, \emph{Geometric structures of vectorial type}, J. Geom. Phys.
56 (2006), 2403-2414.
%
\bibitem{Alexandrov04}
B. Alexandrov, \emph{Hermitian spin surfaces with small eigenvalues of 
the Dolbeault operator}, Ann. Inst. Fourier 54 (2004), 2437-2453.
%
\bibitem{Alexandrov03}
\bysame, \emph{$Sp(n)U(1)$-connections with parallel totally 
skew-symmetric torsion},  Journ. Geom. Phys. 57 (2006), 323-337.
%
\bibitem{Alexandrov&F&S04}
B. Alexandrov, Th. Friedrich, N. Schoemann, \emph{Almost Hermitian 
$6$-manifolds revisited}, J. Geom. Phys. 53 (2005), 1-30.
%
\bibitem{Alexandrov&G&I98}
B. Alexandrov, G. Grantcharov, and S. Ivanov, \emph{An estimate for the
first eigenvalue of the Dirac operator on compact Riemannian spin
manifold admitting parallel one-form}, J. Geom. Phys. 28 (1998), 263-270.
%
%
\bibitem{Alexandrov&G&I01}
B. Alexandrov, G. Grantcharov, and S. Ivanov, \emph{The Dolbeault
operator on Hermitian spin manifolds}, Ann. Inst. Fourier, Univ. Grenoble
51 (2001), 221-225.
%
\bibitem{Alexandrov&I00}
B. Alexandrov, S. Ivanov, \emph{Dirac operators on Hermitian spin surfaces},
Ann.~Global Anal.~Geom. 18 (2000), 529-539.
%
\bibitem{Belgun00}
F. Belgun, \emph{On the metric structure of non-K\"ahler complex surfaces},
Math. Ann. 317 (2000), 1-40.
%
\bibitem{Bismut}
J. M. Bismut, \emph{A local index theorem for non-K\"ahlerian manifolds}, 
Math. Ann. 284 (1989), 681-699.
%
\bibitem{Friedrich80}
Th. Friedrich, \emph{Der erste Eigenwert des Dirac Operators einer 
kompakten Riemannschen Mannigfaltigkeit nichtnegativer Skalarkr\"ummung},
Math. Nachr. 97 (1980), 117-146.
%
\bibitem{Friedrich93}
Th. Friedrich, \emph{The classification of $4$-dimensional K\"ahler 
manifolds with small eigenvalue of the Dirac operator}, Math. Ann. 295 (1993), 565-574.
%
\bibitem{Dirac-Buch}
\bysame, \emph{Dirac operators in Riemannian geometry},
Graduate Studies in Mathematics 25, AMS, Providence, Rhode Island, 2000.
%
\bibitem{Friedrich06}
\bysame, \emph{$\mathrm{G}_2$-manifolds with parallel characteristic torsion},
to appear in Diff. Geom. Appl. 
%
\bibitem{Friedrich&I1}
Th. Friedrich and S. Ivanov, \emph{Parallel spinors and connections with
  skew-symmetric torsion in string theory}, Asian Journ. Math. 6 (2002), 
303-336.
%
\bibitem{Friedrich&I2}
\bysame, \emph{Almost contact manifolds, connections
with torsion and parallel spinors}, J. reine angew. Math. 559 (2003), 
217-236.
%
\bibitem{FriedrichKim}
Th. Friedrich and E.C. Kim, \emph{The Einstein-Dirac equation on 
Riemannian spin manifolds}, J. Geom. Phys. 33 (2000), 128-172.
%
\bibitem{Gauduchon97}
P. Gauduchon, \emph{Hermitian connections and Dirac operators},
Boll. Un. Mat. Ial. ser. VII  2 (1997), 257-289. 
%
\bibitem{Gauduchon&O98}
P. Gauduchon, L. Ornea, \emph{Locally conformally K\"ahler metrics on 
Hopf surfaces}, Ann. Inst. Fourier 48 (1998), 1107-1127.

%
\bibitem{Kirchberg}
K.D. Kirchberg, \emph{An estimate for the first eigenvalue of the Dirac
operator on closed K\"ahler manifolds with positive scalar curvature},
Ann. Glob. Anal. Geom. 4 (1986), 291-326.
%
\bibitem{Kostant99}
B. Kostant, \emph{A cubic {D}irac operator and the emergence of {E}uler number
  multiplets of representations for equal rank subgroups}, Duke Math. J.
 100 (1999), 447-501.
%
%
\bibitem{Schoemann}
N. Schoemann, \emph{Almost Hermitian structures with parallel torsion}, 
to appear.
%
%
\bibitem{Vaisman79}
I. Vaisman, \emph{Locally conformal  K\"ahler manifolds with parallel Lee
form}, Rend. Math. Roma 12 (1979), 263-284.
%
\end{thebibliography}
\end{document}